\documentclass[12pt]{article}
\usepackage{graphicx}
\usepackage{amsfonts}
\usepackage{amssymb}
\usepackage{amsmath}

\font\tendb=msbm10 at 12pt \font\sevendb=msbm10 at 9pt
\font\fivedb=msbm10 at 7pt
\newfam\dbfam
\textfont\dbfam=\tendb \scriptfont\dbfam=\sevendb
\scriptscriptfont\dbfam=\fivedb
\def\db{\fam\dbfam\tendb}
\font\eufm=eufm10\font\eufms=eufm10\font\eufmss=eufm10\newfam\eufam
\textfont\eufam=\eufm\scriptfont\eufam=\eufms\scriptscriptfont\eufam=\eufmss

\font\tendbb=msbm10 at 12pt \font\sevendbb=msbm7 at 9pt
\font\fivedbb=msbm5 at 6pt
\newfam\dbbfam
\textfont\dbbfam=\tendbb \scriptfont\dbbfam=\sevendbb
\scriptscriptfont\dbbfam=\fivedbb

 \def \Z {{\db Z}}

 \def \N {{\db N}}

\def \fin {\hfill \framebox(7,7) }
 
\font\tenMmm=eusm10 at 12pt
\newfam\Mmmfam
\textfont\Mmmfam=\tenMmm

\def\illu #1 by #2 (#3){
  \vbox to #2{
    \hrule width #1 height 0pt depth 0pt
    \vfill
    \special{illustration #3} 
    }
  }

\begin{document}

\null \vspace{2cm}

\begin{center}
{\large {\bf The Tutte  polynomial and the automorphism group of a graph}}\\
 Nafaa Chbili\\
\begin{footnotesize} Department of Mathematical Sciences\\
 College of Science\\
 UAE University\\
 E-mail: nafaachbili@uaeu.ac.ae\\
 \end{footnotesize}

\end{center}
\begin{abstract} A graph $G$  is said to be $p$-periodic, if the  automorphism group $Aut(G)$ contains an element of order $p$  which preserves no edges. In this paper, we investigate  the behavior of graph polynomials (Negmai and  Tutte)   with respect to graph periodicity. In particular,
we prove that if $p$ is a prime, then the coefficients of the Tutte polynomial of such a graph satisfy a certain  necessary condition.
This result is illustrated by an example where the Tutte polynomial  is used to rule out the periodicity of the Frucht graph.
\end{abstract}
\section{Introduction} Throughout this paper, a graph is the geometric realization of a 1-dimensional  CW-complex.
Let $G$ be a graph with edge set $E(G)$ and vertex set $V(G)$. Let $p\geq 2$ be an integer, the graph $G$ is said to be
 \emph{$p$-periodic} if its automorphism group contains an element $h$ of order $p$ subject to the condition  that no edge of $G$ is
fixed under the action of $\Z_p=<h>$. Equivalently, the finite cyclic group $\Z_p$   acts on the graph preserving the incidence
 and such that the restriction of this action to $E(G)$ is fixed point free. Let $G$ be a $p-$periodic graph. If we identify the vertices which belong to the same orbit to a single vertex and the edges which belong to the same orbit to a single edge, then we obtain the  quotient graph  of $G$ under the action of $<h>$. This  quotient graph   is denoted hereafter by $\overline G$.\\

 The chromatic polynomial $P_{G}(\lambda)$ is a classical   tool in graph theory, which counts  the number of proper colorings  of the vertices of the graph with $\lambda$ distinct colors. This polynomial can be recursively defined using a simple contraction-deletion formula. The Tutte polynomial $\tau_{G}(x,y)$  is an isomorphism invariant of graphs. It is a two-variable   polynomial with integral coefficients which specializes to  the chromatic polynomial. The Tutte polynomial can be also  defined through a contraction-deletion formula as explained in Section 2. In this paper, we consider a modified version  of the Tutte polynomial by setting  $T_G(s,t)=\tau_G(s+1,t+1)$. Then $T_G(s,t)=\displaystyle\sum_{i,j}a_{i,j}s^it^j$ where $a_{i,j}$ are integers. The importance of the Tutte polynomial comes not only from  the many
information it carries about the graph, but also from its connection to other fields such as knot theory and  statistical physics. Actually, it is well known that the Tutte polynomial specializes to the partition function of the $q-$state Potts model \cite{WM}.\\
The Tutte polynomial has been generalized into several directions. For instance,  Negami \cite{Ni} introduced  a three variable polynomial $N_G(u,x,y)$  which specializes to the Tutte polynomial (see Section 2). Another interesting generalization  has been obtained by Murasugi \cite{Mu} who defined a polynomial invariant of weighted  graphs, refereed to here as  Murasugi polynomial.  Bollobas and Riordan \cite{BR}  introduced a kind of universal Tutte polynomial of colored graphs  with respect to the deletion-contraction formula.\\
The purpose of this paper is to investigate the relationship  between the automorphism group of a given graph and its polynomials. More precisely, we
study the behavior of the Negami  polynomial of periodic graphs asking whether this polynomial reflects the periodicity of graphs. One of the main motivations of  this study is the nice behavior of the coefficients of the  HOMPLYPT polynomial of
symmetric links introduced in \cite{Ch1} and  \cite{Ch2} for instance. The main results in this paper confirm   that the coefficients of the Negami and Tutte polynomials strongly reflect the periodicity of graphs.

\textbf{Theorem 1.1.} \emph{Let $p$ be a prime and assume  that $G$ is a connected  $p$-periodic graph. Then
$N_G(u,x,y)\cong \displaystyle\sum_{i,j}u^{i}x^{q-jp}y^{jp}$ modulo  $p$, where $q$ is the number of edges of $G$.
}
\\
\textbf{Corollary 1.2.}  \emph{Let $p$ be a prime and assume that $G$ is a connected  $p$-periodic graph. Then
$T(s,t)=\displaystyle\sum_{i,j}a_{i,j}s^it^j$ where $a_{i,j}\cong 0$  modulo $p$  whenever $j-i \ncong 1-r$  modulo $p$, where $r$ is the number of vertices of $G$.\\}

\textbf{ Corollary 1.3.} \emph{Let $G$ be a self-dual connected planar graph and $p$ be  a prime. If  $G$  is $p-$periodic, then $r \cong 1$ modulo $p$, where $r$ is the number of vertices of $G$.\\}

\textbf{Example 1.4.} Let $P$ be the  Petersen graph in Figure 1. It is well  known that this graph is 5-periodic.
The modified Tutte polynomial modulo 5 is:
$T_P(s,t)=s^4 + s^9 + 2 t + 2 s^5 t + s t^2 + t^6$. Since the only nonzero coefficients are $a_{4,0}, a_{9,0}, a_{0,1},a_{5,1},a_{1,2}, a_{0,6}$,  the necessary condition given by the  Theorem is satisfied.
\begin{center}
\includegraphics[width=3cm,height=3cm]{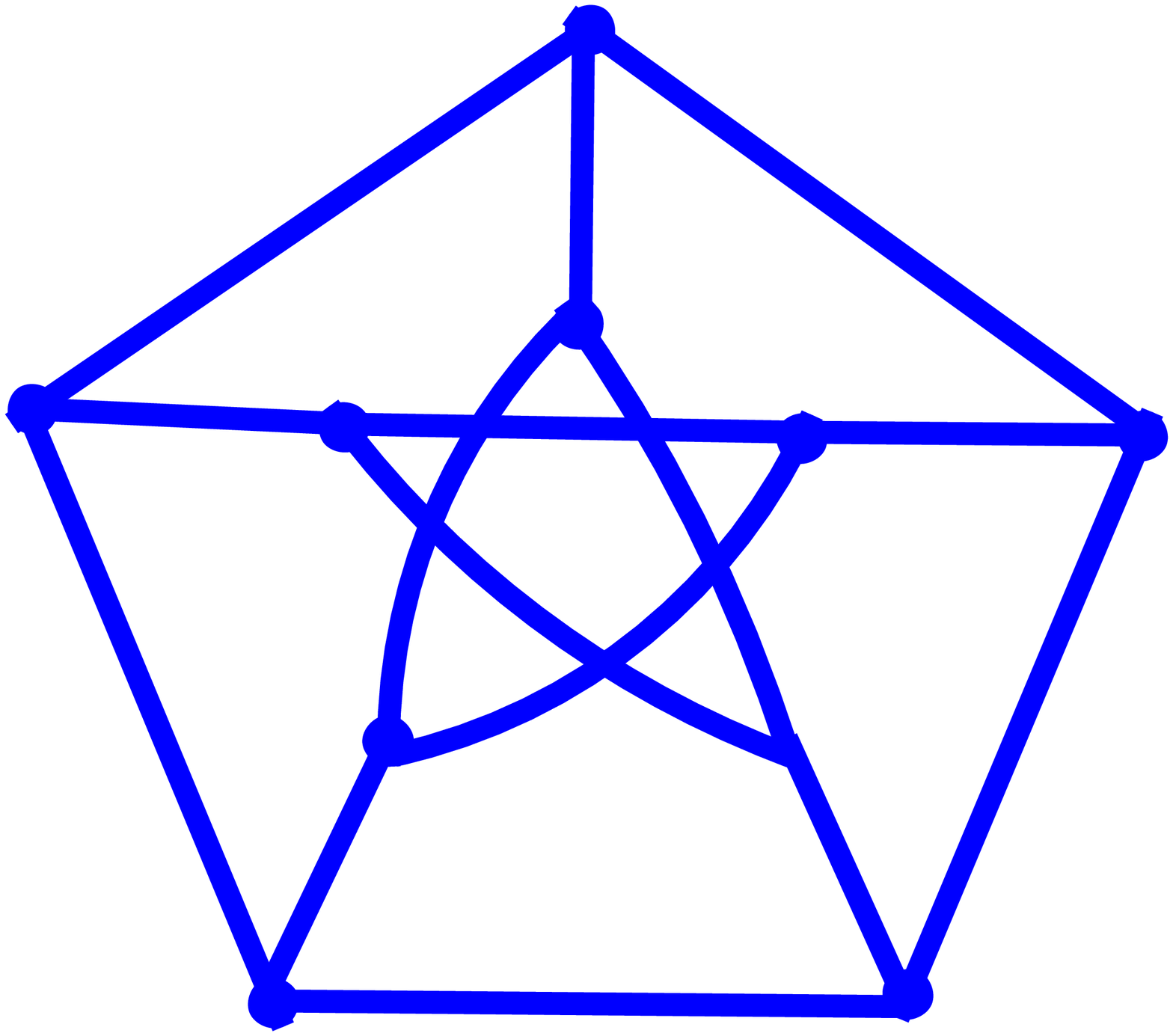}
\end{center}
\begin{center} {\sc  Figure 1} \end{center}
\textbf{Application 1.5.}  Corollary 1.2 can be applied to rule out the possibility of being periodic for certain graphs
 as it is illustrated in the following example.  Let $F$ be the Frucht graph in Figure 2. The reduced modulo 3,  Tutte polynomial of $F$,   is given by \\
$T_F(s,t)=1 + s + s^2 + 2 s^3 + s^6 + 2 s^7 + s^{11} + s t + s^2 t + 2 s^3 t +
 2 s^4 t + s^5 t + 2 s^6 t + 2 s^7 t + s^8 t + 2 t^2 + s^2 t^2 +
 s^4 t^2 + 2 s^5 t^2 + 2 s^6 t^2 + s^7 t^2 + t^3 + 2 s t^3 + s^2 t^3 +
  s^4 t^3 + s^5 t^3 + 2 t^4 + 2 s t^4 + 2 s^2 t^4 + 2 s^3 t^4 +
 s t^5 + t^7.$
Since $T_F(s,t)$ does not satisfy the periodicity condition given by Corollary 1.2.  Then the  Frucht graph is not $3-$periodic.
\begin{center}
\includegraphics[width=3cm,height=3cm]{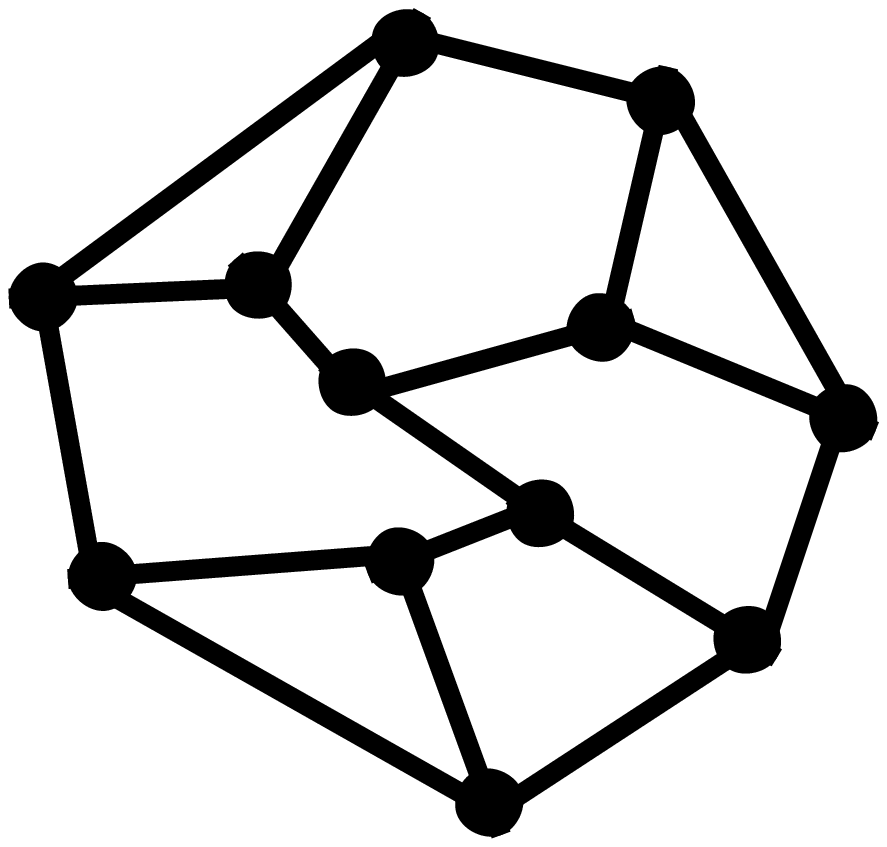}
\end{center}
\begin{center} {\sc  Figure 2} \end{center}

\section{Negami's polynomial}
The Negami polynomial is a three-variable isomorphism invariant of graphs \cite{Ni}. It can be defined by a simple
expansion formula. Let $G$ be a graph with edges set $E(G)$. Assume that $G$ has $r$ vertices, $q$ edges and $w(G)$ components.
 Then the Negami polynomial is defined by the following formula:
$$ N_G(u,x,y)=\displaystyle\sum_{Y\subset E(G)}u^{w(G-Y)}x^{q-|Y|}y^{|Y|}.$$

According to \cite{Ni},  the Tutte polynomial is obtained as a specialization of the Negami polynomial:
$$\tau_G(x,y)=(y-1)^{-r}(x-1)^{-w(G)}N((x-1)(y-1),1,y-1).$$ Thus we have
$$T_G(s,t)=\tau(s+1,t+1)=t^{-r}s^{-w(G)}N(st,1,t).$$

If $G$ is a graph and $e$ is an edge, we denote by $G-e$ the graph obtained from $G$ by deleting the edge $e$ and by $G/e$ we denote the graph obtained by contracting the edge $e$. The Tutte polynomial can be  defined using the following  deletion-contraction formula:

$$ \tau_G(x,y) = \left\{ 
\begin{array}{llll}
x\tau_{G/e}(x,y)&&& \mbox { if $e$ is a bridge} \\
 y\tau_{G-e}(x,y) &&& \mbox {if $e$ is a loop}\\
\tau_{G-e}(x,y)+\tau_{G/e}(x,y) &&& \mbox{if $e$ is an ordinary edge}
\end{array} \right .
 $$

 Together with the condition $T(E_n)=1$, where $E_n$ is the graph with $n$ vertices and no edges.

\textbf{Proof of Theorem 1.1.} Let $G$ be a $p-$periodic connected graph with $r$ vertices and  $q$ edges.
The cyclic group $\Z_p$ acts on ${\mathcal P} (E(G))$. Since $p$ is a prime then there are two kinds of orbits.
Orbits which are made up of only one element and orbits which are made up of exactly $p$ elements.\\
It is clear that if $Y$ belongs to an orbit made up of $p$ elements, then the contribution of these elements
sums to zero modulo $p$  in the expansion formula of $N_G(u,x,y)$. Consequently, in the computation of $N_G(u,x,y)$ modulo $p$, we shall
 consider only subsets $Y$ which are fixed by the action of $\Z_p$. Since the action has no fixed edges, then the number of elements of $|Y|$ is an
multiple of $p$, say $jp$. A priori, we don't have any control of the number of components  $w(G-Y)$. Thus, the the polynomial is written on the form:\\
$$ N_G(u,x,y)=\displaystyle\sum_{Y\subset E(G)}u^{w(G-Y)}x^{q-|Y|}y^{|Y|} \cong
\displaystyle\sum_{i,j\in S} u^{i}x^{q-jp}y^{jp} \mbox{ modulo } p,$$
where $S$ is a certain subset of $\N\times\N$. \fin \\

\textbf{Proof of Corollary 1.2.} Since $G$ is a connected graph, then  $T_G(s,t)=t^{-r}s^{-1}N_G(st,1,t)$.
 The congruence condition in Theorem 1.1 implies that modulo $p$, we have  $T_G(s,t)\cong t^{-r}s^{-1}\displaystyle\sum_{i,j}(st)^it^{jp}$. Hence, we get
$T_G(s,t)\ \displaystyle\sum_{i,j}s^{i-1}t^{jp-r+i}$.
 This ends the proof of Corollary 1.2. \fin \\

\textbf{Proof of Corollary 1.3.} Let $G^{*}$ be the dual graph of $G$. It is well known that $T_{G}(s,t)=T_{G^{*}}(t,s)$. Thus, if  $G$ is
self-dual then $T_{G}(s,t)=T_{G}(t,s)$. If we assume that $G$ is $p-$periodic then $T(s,t)=\displaystyle\sum_{i,j}a_{i,j}s^it^j$ where\\
$a_{i,j}\cong 0$  modulo $p$  whenever $j-i \ncong 1-r$  modulo $p$ and \\
$a_{i,j}\cong 0$  modulo $p$  whenever $i-j \ncong 1-r$  modulo $p$.\\
This implies  that $a_{i,j}\cong 0$  modulo $p$  whenever $r\ncong 1 $  modulo $p$. \fin \\

\section{Negami's polynomial and Murasugi congruence}
In Knot theory, Murasugi's congruence refers to the relationship between the invariant of a periodic link and the invariant of its quotient link.
 Murasugi proved his congruence for the Alexander polynomial of periodic knots \cite{Mu1}, then extended  it to the Jones polynomial
 \cite{Mu} of periodic links. This congruence has been generalized to other link invariants and to the Yamada invariant
of  spatial graphs, \cite{Ch3,Ch4}. In \cite{Mu}, Murasugi introduced a new polynomial of weighted graphs and proved that
this polynomial satisfies a certain congruence of Murasugi type. In this paragraph, we prove that this condition extends to
the Negami and Tutte polynomials as well. \\

\textbf{Theorem 3.1.} \emph{Let $p$ be a prime and $G$ be a connected $p-$periodic graph. Then $N_G(u,x,y) \cong (N_{\bar G}(u,x,y))^p$ modulo the ideal generated
by $p$ and  $u^p-u$.}\\
\textbf{Corollary 3.2.} \emph{Let $p$ be a prime and $G$ be a connected $p-$periodic graph. Then, $T_G(s,t) \cong (T_{\bar G}(s,t))^p$ modulo
the ideal generated
by $p,$ $s^p-s$ and  $t^p-t$.}\\

\textbf{Proofs.} A subset $Y \in {\mathcal P}(E(G))$ is called periodic if and only if $Y$ is
invariant by the $\Z_p$-action on $E(G)$. As we have seen in the proof of
Theorem 1.1, only periodic subsets $Y \in {\mathcal P}(E(G))$ are to  be considered in the computation
of $N_G(u,x,y)$ modulo $p$. It can be easily seen that a periodic subset  $Y$  defines a subset $\overline Y$ of $E(\overline G)$ and vise-versa.
 Moreover,
we have $|Y|=p|\overline Y|$ and  $q=p\overline q$ where $\overline q$ is the number of edges of the quotient graph $\overline G$.\\
On the other hand, it can be easily seen that if  $\overline G-\overline Y$ has one component, then $G-Y$ can have one or $p$ components. Hence,
$u^{w(G-Y)}$ and $u^{w(\overline G-\overline Y)}$ coincide modulo $u^p-u$. Now, modulo $p$ and $u^p-u$, the Negami polynomial of $G$ can be
 computed as follows:
$$N_G(u,x,y)\cong \displaystyle\sum_{Y\subset E(G)}u^{w(G-Y)}x^{q-|Y|}y^{|Y|}\cong
\displaystyle\sum_{\overline Y \subset E(\overline G)}u^{w(\overline G-\overline Y)}x^{p\overline q-p|\overline Y|}y^{p|\overline Y|}\cong
[N_{\overline G}(u,x,y)]^p.$$

The proof of Corollary 3.2 is straightforward. \fin \\

\textbf{Remark.} Since the periodicity criterion  introduced in  Theorem 3.1 involves both the  polynomials of the graph and its quotient. Then it
seems difficult  to find an example where this condition is  used to rule out the periodicity of the graph.
However, this might be possible if we add some extra conditions on the graph or its quotient as in the following case.
 Assume that $G$ is a $p$-periodic graph having a cycle of order $p$ which is invariant by the action of $\Z_p$.
Let us first write the criterion of Theorem 3.1 for the  special case of the chromatic polynomial. According to \cite{Ni},  we have $P_{G}(\lambda)=N_G(\lambda,-1,1)$. Then  Theorem 3.1  implies
$$P_{G}(\lambda) \cong (P_{\overline  G}(\lambda))^p \;\;\mbox{  modulo the ideal generated
by } p \mbox{ and } \lambda^p-\lambda.$$ \\
If  $G$ has a cycle of length  $p$ which is invariant under the action of $\Z_p$ then the quotient graph has a loop and thus
 $P_{\bar G}(\lambda)=0$. Hence: $$P_{G}(\lambda) \cong 0\;\;\mbox{  modulo the ideal generated
by } p \mbox{ and } \lambda^p-\lambda.$$

\end{document}